\documentclass[11pt,a4paper]{article}
\usepackage{a4wide}
\usepackage{amsmath,amssymb}
\usepackage[english]{babel}
\usepackage{graphicx}
\usepackage[numbers]{natbib}
\usepackage{theorem}
\usepackage{bm}
\usepackage{titlesec}
\usepackage{ifpdf}
\usepackage{hyperref}
\usepackage{setspace}

\hypersetup{
    bookmarksopen=false,
    bookmarksnumbered=true,
    pdftitle={Closed-Form Waiting Time Approximations for Polling Systems},
    pdfauthor={M.A.A. Boon, E.M.M. Winands, I.J.B.F. Adan, A.C.C. van Wijk}
}
\ifpdf
  \hypersetup{colorlinks=true,linkcolor=black,urlcolor=black,citecolor=black,pdfpagemode=UseOutlines,plainpages=false,pdfpagelabels}
\else
  \hypersetup{colorlinks=false}
\fi

\titleformat{\subsubsection}[hang]{\normalfont\em}{}{1ex}{}

\theorembodyfont{\upshape}
\theoremheaderfont{\bfseries}

\newtheorem{theorem}{Theorem}[section]

\newtheorem{remark}[theorem]{Remark}
\newtheorem{example}{Example}

\numberwithin{equation}{section}

\newcommand{\ee}{\textrm{e}}
\newcommand{\dd}{\,\textrm{d}}

\newcommand{\E}{\mathbb{E}}
\newcommand{\V}{\mathbb{V}\textrm{ar}}

\renewcommand{\O}{\mathcal{O}}
\newcommand{\C}{\textit{cv}}

\providecommand{\href}[2]{#2}

\title{Closed-Form Waiting Time Approximations for Polling Systems\footnote{The research was done in the framework of the BSIK/BRICKS project, and of the European Network of Excellence Euro-NF.}}
\author{M.A.A. Boon\footnote{\textsc{Eurandom} and Department of Mathematics and Computer Science, Eindhoven University of Technology, P.O. Box 513, 5600MB Eindhoven, The Netherlands}\\\href{mailto:marko@win.tue.nl}{marko@win.tue.nl} \and E.M.M. Winands \footnote{Department of Mathematics, Section Stochastics, VU University, De Boelelaan 1081a, 1081HV Amsterdam, The Netherlands}\\\href{mailto:emm.winands@few.vu.nl}{emm.winands@few.vu.nl} \and I.J.B.F. Adan\footnotemark[2]\\\href{mailto:iadan@win.tue.nl}{iadan@win.tue.nl} \and A.C.C. van Wijk\footnote{Department of Industrial Engineering \& Innovation Sciences and Department of Mathematics and Computer Science, Eindhoven University of Technology, P.O. Box 513, 5600MB Eindhoven, The Netherlands} \\\href{mailto:a.c.c.v.wijk@tue.nl}{a.c.c.v.wijk@tue.nl}}
\date{December, 2010}

\begin{document}
\maketitle

\begin{abstract}
A typical polling system consists of a number of queues, attended by a single server in a fixed order.
The vast majority of papers on polling systems focusses on Poisson arrivals, whereas very few results are available for general arrivals.
The current study is the first one presenting simple closed-form approximations for the mean waiting times in polling systems with renewal arrival processes, performing well for \emph{all} workloads. The approximations are constructed using heavy traffic limits and newly developed light traffic limits. The closed-form approximations may prove to be extremely useful for system design and optimisation in application areas as diverse as telecommunication, maintenance, manufacturing and transportation.

\end{abstract}

\section{Introduction}\label{introduction}

Polling systems are queueing systems consisting of multiple queues, visited by a single server - typically in a fixed, cyclic order. They find their origin in many real-life applications, e.g. (computer) communication, production and manufacturing environments, traffic and transportation. For a good literature overview of polling systems and their applications, we refer to surveys of, e.g., Takagi \cite{takagi1988qap}, Levy and Sidi \cite{levysidi90}, and Vishnevskii and Semenova \cite{vishnevskiisemenova06}. When studying literature on polling systems, it rapidly becomes apparent that the computation of the distributions and moments of the waiting times and marginal queue lengths is very cumbersome. Closed form expressions do not exist, and even when one specifies the number of queues and solves the set of equations that leads to the mean waiting times, the obtained expressions are still too lengthy and complicated to interpret directly.
Numerical procedures, both approximate and exact, have been developed in the past to compute these performance measures. However, these methods have several drawbacks. Firstly, they are not transparent and act as a kind of black box. It is, for instance, rather difficult to study the impact of parameters like the occupation rate and the service level. Secondly, these procedures are computationally complex and hard, if not impossible, to implement in a standard spreadsheet program commonly used on the work floor. Finally, the vast majority of standard methods focusses on Poisson arrival processes, which may not be very realistic in many application areas.
In the present paper we study polling systems in which the arrival streams are not (necessarily) Poisson, i.e., the interarrival times follow a general distribution.
The goal is to derive closed-form approximate solutions for the mean waiting times and mean marginal queue lengths, which can be computed by simple spreadsheet calculations. 

Our approach in developing an approximation for the mean waiting times uses novel developments in polling literature. Recently, a heavy traffic (HT) limit has been developed for the mean waiting times as the system becomes saturated \cite{coffman98,olsenvdmei05,vdmeiwinands08}. In the present paper we derive an approximation for the light traffic (LT) limit, i.e. as the load decreases to zero, which is exact for Poisson arrivals. The main idea is to create an interpolation between the LT limit and the HT limit. This interpolation yields good results, and has several nice properties, like satisfying the Pseudo Conservation Law (PCL), and being exact for symmetric systems with Poisson arrivals and in many limiting cases. These properties are described in more detail in the present paper. In polling literature, several alternative approximations have been developed before, most of which assume Poisson arrivals.
For polling systems with Poisson arrivals and gated or exhaustive service, the best results, by far, are obtained by an approximation based on the PCL  (see, e.g., \cite{boxmaworkloadsandwaitingtimes89,everitt86,groenendijk89}). Fischer et al. \cite{fischer2000} study an approximation for the mean waiting times in polling systems, which is also based on an interpolation between (approximate) LT and HT limits. Their approach, however, is applied to a system with Poisson arrivals and time-limited service.
Hardly any closed-form approximations exist for non-Poisson arrivals. The few that exist, perform well in specific limiting cases, e.g., under HT conditions \cite{olsenvdmei05,vdmeiwinands08}, or if switch-over times become very large \cite{winandslargesetups07,winandslargesetupsbranching09}, but performance deteriorates rapidly if these limiting conditions are abandoned, in contrast to the approximation developed in the present paper. We show in an extensive numerical study that the quality of our approximation can be compared to the PCL approximation for systems with Poisson arrivals, but provides good results as well for systems with renewal arrivals.

Because of its simple form, the approximation function is very suitable for optimisation purposes and implementation in a spreadsheet. Although only the mean waiting times of systems with exhaustive or gated service are studied, the results can be extended to higher moments and general branching-type service disciplines. Polling systems with polling tables and/or batch service can also be analysed in a similar manner.

The structure of the present paper is as follows: the next section introduces the model and the required notation, and states the main result. Section \ref{approximationsection} illustrates how this main result is obtained, while Section \ref{numericalexample} provides results on the accuracy of the approximation for a large set of combinations of input parameter values. The last section discusses further research topics and possible extensions of the model.

\section{Model description and main result}\label{modeldescription}

The model under consideration is a polling system consisting of $N$ queues, $Q_1, \dots, Q_N$, with renewal arrival processes. Indices throughout the present paper are understood to be modulo $N$: $Q_{N+1}$ actually refers to $Q_1$. Whenever a server switches from $Q_i$ to $Q_{i+1}$, a random switch-over time $S_i$ is incurred. The generic service requirement of a customer arriving in $Q_i$, also referred to as a type $i$ customer, is denoted by the random variable $B_i$.
We make the usual independence assumptions for polling systems; the interarrival times, service times and switch-over times are all independent.
The moment at which the server switches from one queue to the next queue, is determined by the \emph{service discipline} of the queue that is being served. In the present paper we focus on polling systems in which each queue is either served according to the \emph{gated} service discipline, which states that during the course of a visit of the server to $Q_i$, only those type $i$ customers are served that were present at the beginning of that visit, or according to the \emph{exhaustive} service discipline, which means that the server keeps on serving type $i$ customers until $Q_i$ is empty, before switching to $Q_{i+1}$.

We regard several variables as a function of the load $\rho$ in the system. Scaling is done by keeping the service time distributions fixed, and varying the interarrival times. For each variable $x$ that is a function of the load in the system, $\rho$, its value evaluated at $\rho=1$ is denoted by $\hat x$. For $\rho = 1$, the generic interarrival time of the stream in $Q_i$ is denoted by $\hat{A}_i$. Reducing the load $\rho$ is done by scaling the interarrival times, i.e., taking the random variable $A_i := \hat{A}_i/\rho$ as generic interarrival time at $Q_i$. After scaling, the load at $Q_i$ becomes $\rho_i = \rho\frac{\E[B_i]}{\E[\hat{A}_i]}$. The (scaled) rate of the arrival stream at $Q_i$ is defined as $\lambda_i = 1/\E[A_i]$. Similarly, we define arrival rates $\hat\lambda_i = 1/\E[\hat{A}_i]$, and proportional load at $Q_i$, $\hat\rho_i = \frac{\rho_i}{\rho}$ (``proportional'' because $\sum_{i=1}^N \hat\rho_i = 1$). The system is assumed to be stable, so $\rho$ is varied between 0 and 1.

We use $B$ to denote the generic service requirement of an arbitrary customer entering the system, with $\E[B^k] = \frac{\sum_{i=1}^N \hat\lambda_i \E[B_i^k]}{\sum_{j=1}^N \hat\lambda_j}$ for any integer $k > 0$, and $S = \sum_{i=1}^N S_i$ denotes the total switch-over time in a cycle. Finally, the (equilibrium) residual length of a random variable $X$ is denoted by $X^\textit{res}$, with $\E[X^\textit{res}] = \frac12\E[X^2]/\E[X]$.

We now present the main result of this paper, which is a closed-form approximation formula for the mean waiting time $\E[W_i]$ of a type $i$ customer as a function of $\rho$:
\begin{equation}
\E[W_{i,\textit{app}}]=\frac{K_{0,i}+K_{1,i} \rho + K_{2,i} \rho^2}{1-\rho},\qquad i=1,\dots,N.\label{ewapprox}
\end{equation}
The constants $K_{0,i}, K_{1,i}$, and $K_{2,i}$ depend on the input parameters and the service discipline. If all queues receive \emph{exhaustive} service, the constants become:
\begin{align}
K_{0,i} =& \E[S^{\textit{res}}],\label{c0}\\
K_{1,i} =& \hat\rho_i \big(\E[\hat{A}_i]\hat{g}_i(0) -1\big)\E[B_i^\textit{res}] + \E[B^{\textit{res}}] +\hat\rho_i\big(\E[S^{\textit{res}}]-\E[S]\big) -\frac{1}{\E[S]}\sum_{j=0}^{N-1}\sum_{k=0}^j\hat\rho_{i+k}\V[S_{i+j}],\label{c1}\\
K_{2,i} =& \frac{1-\hat{\rho}_i}{2} \left( \frac{\sum_{j=1}^N \hat{\lambda}_j \left( \V[B_j] + \hat{\rho}_j^2 \V[\hat{A}_j] \right)}{\sum_{j=1}^N \hat{\rho}_j (1-\hat{\rho}_j)} +  \E[S] \right) - K_{0,i} - K_{1,i}.\label{c2}
\end{align}
If all queues receive \emph{gated} service, we get:
\begin{align}
K_{0,i} =& \E[S^{\textit{res}}],\label{c0gated}\\
K_{1,i} =& \hat\rho_i \big(\E[\hat{A}_i]\hat{g}_i(0) -1\big)\E[B_i^\textit{res}] + \E[B^{\textit{res}}] +\hat\rho_i\E[S^{\textit{res}}] -\frac{1}{\E[S]}\sum_{j=0}^{N-1}\sum_{k=0}^j\hat\rho_{i+k}\V[S_{i+j}],\hphantom{\big(-\E[S]\big)}\label{c1gated}\\
K_{2,i} =& \frac{1+\hat{\rho}_i}{2} \left( \frac{\sum_{j=1}^N \hat{\lambda}_j \left( \V[B_j] + \hat{\rho}_j^2 \V[\hat{A}_j] \right)}{\sum_{j=1}^N \hat{\rho}_j (1+\hat{\rho}_j)} +  \E[S] \right) - K_{0,i} - K_{1,i}.\label{c2gated}
\end{align}
The term $\hat{g}_i(t)$ is the density of $\hat{A}_i$, the interarrival times at $\rho=1$. This term is discussed in more detail in the next section, but for practical purposes it is useful to know that $\E[\hat{A}_i]\hat{g}_i(0)$ can be very well approximated by
\[
\E[\hat{A}_i]\hat{g}_i(0) \approx \begin{cases}
2\frac{\C_{A_i}^2}{\C_{A_i}^2+1} &\qquad\text{ if } \C_{A_i}^2>1,\\
\left(\C_{A_i}^2\right)^4 &\qquad\text{ if } \C_{A_i}^2\leq1,
\end{cases}
\]
where $\C_{A_i}^2$ is the squared coefficient of variation (SCV) of $A_i$ (and, hence, also of $\hat{A}_i$). Note that this simplification results in an approximation that requires only the first two moments of each input variable (i.e., service times, switch-over times, and interarrival times).

\begin{remark}
In case of Poisson arrivals, the constants $K_{1,i}$ and $K_{2,i}$ simplify considerably. E.g., for exhaustive service they simplify to:
\begin{align*}
K_{1,i}^\textit{Poisson} =& \E[B^{\textit{res}}] +\hat\rho_i\big(\E[S^{\textit{res}}]-\E[S]\big) -\frac{1}{\E[S]}\sum_{j=0}^{N-1}\sum_{k=0}^j\hat\rho_{i+k}\V[S_{i+j}],\\
K_{2,i}^\textit{Poisson} =& (1-\hat{\rho}_i) \left( \frac{\E[B^\textit{res}]}{\sum_{j=1}^N \hat{\rho}_j (1-\hat{\rho}_j)} +  \frac{\E[S]}{2} \right) - K_{0,i} - K_{1,i}^\textit{Poisson}.
\end{align*}
\end{remark}

The derivation of this approximative formula for the mean waiting time is the topic of the next section. An approximation for the mean \emph{queue length} at $Q_i$, $\E[L_i]$ is obtained by application of Little's Law to the \emph{sojourn time} of type $i$ customers, i.e. the waiting time plus the service time. As a function of $\rho$, we have
\[\E[L_{i,\textit{app}}] = \rho\frac{\E[W_{i,\textit{app}}]+\E[B_i]}{\E[\hat{A}_i]}.\]

\section{Derivation of the approximation}\label{approximationsection}

Approximation (\ref{ewapprox}) is an interpolation approximation based on LT and HT limits. In the next subsection we first provide a motivation for this approach.

\subsection{Generic interpolation function}

In its generic form, the interpolation approximation proceeds as follows; see \cite{reimansimon88, simon92}. Consider an open queueing system with load $\rho$. Let $f(\rho)$, $0 \le \rho < 1$, be some function of the queueing system (such as the mean waiting time), which is assumed to be analytic on $[0, 1)$, i.e., by Taylor's Theorem $f(\rho)$ can be expressed as
\[
f(\rho) = \sum_{n=0}^\infty \frac{f^{(n)}(0)}{n!} \rho^n , \quad 0 \le \rho < 1,
\]
where $f^{(n)}(\rho)$ denotes the $n$th derivative $f(\rho)$. Usually, $f(\rho)$ is intractable, but it may be possible to derive partial information about $f(\rho)$, such as the light traffic limits $f^{(n)} (0)$ for $n=0, 1, \ldots, k$ and the ``canonical'' heavy traffic limit
$h = \lim_{\rho \rightarrow 1} (1 - \rho)f(\rho)$. For examples, see \cite{fleming1991, reimansimon88, whitt89}, where based on the partial information, an approximation for $f(\rho)$ is constructed of the form
\begin{equation}\label{inter1}
\tilde{f}(\rho) = \frac{q(\rho)}{1-\rho} ,
\end{equation}
where $q(\rho)$ is the $(k+1)$st degree polynomial, uniquely determined by the requirement that $\tilde{f}(\rho)$ has to match everything that is known about $f(\rho)$, i.e., $\tilde{f}^{(n)}(0) = f^{(n)}(0)$ for $n = 0, 1, \ldots, k$ and
$\lim_{\rho \rightarrow 1} (1 - \rho)\tilde{f}(\rho) = h$.
The heavy traffic limit implies (see Prop. 1 in \cite{simon92}),
\[
\lim_{n \rightarrow \infty} \frac{f^{(n)}(0)}{n!} = h .
\]
This suggests that, in the Taylor series, $f^{(n)}(0)$ for $n > k$ can be approximated by $n! h$. Thus, combined with knowledge of the light traffic limits $f^{(n)} (0)$ for $n=0, 1, \ldots, k$, the following new approximation can be produced,
\[
\bar{f} (\rho) = \sum_{n=0}^k \frac{f^{(n)}(0)}{n!} \rho^n + h \frac{\rho^{k+1}}{1-\rho} .
\]
Interestingly, this seemingly different approximation $\bar{f}(\rho)$ is identical to the interpolation approximation $\tilde{f}(\rho)$, confirming the notion that they are the ``natural''
approximation for $f(\rho)$, given the partial information. In \cite{reimansimon88} the interpolation approximation (\ref{inter1}) is shown to work extremely well for several examples.
The present paper applies approximation (\ref{inter1}) to the new setting of polling systems with general renewal arrivals, for which no analytic expressions are known for the mean waiting times. Choosing $f(\rho)$ as the mean waiting time of a type $i$ customer, we derive new (approximations for the) light traffic limits $f(0)$ and the first derivative $f'(0)$, which together with the heavy traffic limit, yield an interpolation with a quadratic polynomial $q(\rho)$; see (\ref{ewapprox}).

In Sections \ref{lighttraffic} and \ref{heavytraffic} we derive the LT and HT limits, respectively. The interpolation approximation, matching these limits, is presented in Section \ref{inter}. Finally, in Section \ref{specialcases}, we show that the interpolation approximation also matches known exact results for mean waiting times in polling systems.
In particular, we prove that the form \eqref{inter1} is crucial for satisfying the \emph{pseudo-conservation law} for all loads $\rho$.
An extensive numerical validation is the topic of Section \ref{numericalexample}, showing that the interpolation approximation works extremely well for polling systems.

\subsection{Light traffic}\label{lighttraffic}

The mean waiting times in the polling model under consideration in light-traffic, have been studied in Blanc and Van der Mei \cite{blancvdmei95}, under the assumption of Poisson arrivals. They obtain expressions for the mean waiting times in light traffic that are exact up to (and including) first-order terms in $\rho$. These expressions have been found by carefully inspecting numerical results obtained with the Power-Series Algorithm, but no proof is provided. In the present section we shall not only prove the correctness of the light-traffic results in a system with Poisson arrivals, but also use them as base for an approximation for the mean waiting times in polling systems with renewal interarrival times.
The key ingredient to the LT analysis of a polling system, is the well-known Fuhrmann-Cooper decomposition \cite{fuhrmanncooper85}. It states that in a vacation system with Poisson arrivals the queue length of a customer is the sum of two independent random variables: the number of customers in an isolated $M/G/1$ queue, and the number of customers during an arbitrary moment in the vacation period. The distributional form of Little's Law \cite{keilsonservi90} can be used to translate this result to waiting times. Since no independence is required between the length of a vacation and the length of the preceding visit period, this decomposition also holds for polling systems with Poisson arrivals. We introduce $V_i$ to denote the length of a visit period to $Q_i$, and $I_i$ to denote the length of the intervisit period, i.e. the time that the server is away between two successive visits to $Q_i$. Using $C_i$ to denote the cycle time, starting at a visit beginning to $Q_i$, we have $\E[V_i] = \rho_i\E[C_i]$ and $\E[I_i] = (1-\rho_i)\E[C_i]$. It is well-known that the mean cycle time in polling systems, unlike higher moments, does not depend on the starting point: $\E[C_i]=\E[C]=\frac{\E[S]}{1-\rho}$.

The Fuhrmann-Cooper decomposition, applied to the mean waiting time, results in:
\begin{align}
&\textrm{exhaustive: } & \E[W_i] &= \E[W_{i,M/G/1}]+\E[I_i^\textit{res}],\label{decompositionexhaustive}\\
&\textrm{gated: } & \E[W_i] &= \E[W_{i,M/G/1}]+\E[I_i^\textit{res}]+\frac{\E[V_iI_i]}{\E[I_i]}.\label{decompositiongated}
\end{align}
For our approximation, we assume that this decomposition also holds for renewal arrival processes in light traffic. Determining the LT limit of the mean waiting time, $\E[W_i^{\textit{LT}}]$,  in a polling system with exhaustive or gated service is based on the following two-step approach. The first step is to find the LT limit of $\E[W_{i,GI/G/1}]$, the mean waiting time of a $GI/G/1$ queue with only type $i$ customers in isolation, $i=1,\dots,N$. The second step is determining $\E[I_i^{\textit{res}}]$, the mean residual intervisit time of $Q_i$, and $\frac{\E[V_iI_i]}{\E[I_i]}$, the mean visit time of $Q_i$ given that it is being observed at a random epoch during the following intervisit time. 

\begin{remark}
Bertsimas and Mourtzinou \cite{bertsimas99} state that the decompositions \eqref{decompositionexhaustive} and \eqref{decompositiongated} also hold for polling systems with Mixed Generalised Erlang arrivals. However, simulation and exact analysis of some simple cases indicate that the decomposition result is not valid for the mean waiting times.
\end{remark}

For the LT limit of the mean waiting time in a $GI/G/1$ queue, we use Whitt's result (Equation (16) in \cite{whitt89}), which gives:
\begin{equation}
\lim_{\rho_i\downarrow0} \frac{\E[W_{i,GI/G/1}]}{\rho_i}= \frac{1+\C_{B_i}^2}{2}\E[\hat{A}_i]\hat{g}_i(0)\E[B_i],\label{EWwhitt}
\end{equation}
where $\C_{B_i}^2$ is the SCV of the service times, and $\hat{g}_i(t)$ is the density of the interarrival times $\hat{A}_i$.
For practical purposes, it may be more convenient to express $\hat{g}_i(0)$ in terms of the density of $A_i$, the generic interarrival time of $Q_i$ in the scaled situation.
The relation between the density of the scaled interarrival times $A_i\, (= \hat{A}_i/\rho)$, denoted by $g_{i}(t)$, and the density of $\hat{A}_i$, $\hat{g}_i(t)$, is simply: $g_{i}(t) = \rho \hat{g}_i(\rho t)$. This means that the term $\E[\hat{A}_i]\hat{g}_i(0)$ can be rewritten as
\[\E[\hat{A}_i]\hat{g}_i(0)=\E[A_i]g_{i}(0).\]
Because of this equality, in the remainder of the paper we might use either notation. Since determining $\E[\hat{A}_i]\hat{g}_i(0)$ is a required step in the computation of our approximation for $\E[W_i]$, we give some practical examples.
\begin{example}\label{exampleC1}
If the scaled interarrival times $A_i$ are exponentially distributed with parameter $\lambda_i := 1/\E[A_i]$, we have $g_{i}(t) = \lambda_i \ee^{-\lambda_i t}$. This implies that $\E[A_i]g_i(0) = 1$.
\end{example}
\begin{example}\label{exampleCgeq1}
In this example we assume that $A_i$ follows a $H_2$ distribution with balanced means. The SCV of $A_i$ is denoted by $\C_{A_i}^2$. The density of this hyper-exponential distribution is (see, e.g., \cite{tijms94})
\[g_{i}(t) = p \mu_1 \ee^{-\mu_1 t} + (1-p) \mu_2 \ee^{-\mu_2 t},\]
with
\begin{eqnarray*}
p &=& \frac12\left(1 + \sqrt{\frac{\C_{A_i}^2-1}{\C_{A_i}^2+1}}\right),\\
\mu_1 &=& \frac{1}{\E\left[A_i\right]}\left(1 + \sqrt{\frac{\C_{A_i}^2-1}{\C_{A_i}^2+1}}\right),\\
\mu_2 &=& \frac{1}{\E\left[A_i\right]}\left(1 - \sqrt{\frac{\C_{A_i}^2-1}{\C_{A_i}^2+1}}\right).\\
\end{eqnarray*}
This leads to $\E[A_i]g_i(0) =  1 + \frac{\C^2_A-1}{\C^2_A+1} = 2\frac{\C^2_A}{\C^2_A+1}$. 
\end{example}
\begin{example}\label{exampleCleq1}
Now we assume that the interarrival times follow a mixed Erlang distribution. The density of the scaled interarrival times is:
\[g_{i}(t) = p \frac{\mu^{k-1}t^{k-2}}{(k-2)!}\ee^{-\mu t} + (1-p) \frac{\mu^kt^{k-1}}{(k-1)!}\ee^{-\mu t},\]
i.e., a mixture of an Erlang($k-1$) and an Erlang($k$) distribution with
\begin{eqnarray*}
k &=& \left\lceil\frac{1}{\C_{A_i}^2}\right\rceil, \\
p &=& \frac{k\,\C_{A_i}^2-\sqrt{k(1+\C_{A_i}^2)-k^2\,\C_{A_i}^2}}{1+\C_{A_i}^2},\\
\mu &=& \frac{k-p}{\E[A_i]}.
\end{eqnarray*}
If $k>2$, this leads to $\E[A_i]\,g_{i}(0) =  0$.  
\end{example}
The distributions in Examples $\ref{exampleC1}-\ref{exampleCleq1}$ are typical distributions to be used in a two-moment fit if the SCV of the interarrival times is respectively 1, greater than 1, and less than 1  (cf. \cite{tijms94}).
The examples illustrate how $\E[A_i]g_i(0)$ can be computed if the density of the (scaled) interarrival times is known. If no information is available about the
complete density, but the first two moments of $A_i$ are known, Whitt suggests to use the following approximation for $\E[A_i]g_i(0)$:
\[
\E[A_i]g_i(0) = \begin{cases}
2\frac{\C_{A_i}^2}{\C_{A_i}^2+1} \qquad&\text{ if } \C_{A_i}^2>1, \\
\left(\C_{A_i}^2\right)^4 &\text{ if } \C_{A_i}^2\leq1,
\end{cases}
\]
where $\C_{A_i}^2$ is the squared coefficient of variation of the interarrival times of $Q_i$. This approximation is exact for $\C_{A_i}^2>1$, if the interarrival time distribution is a hyper-exponential distribution as discussed in Example \ref{exampleCgeq1}. For $\C_{A_i}^2\leq1$, the approximation is rather arbitrary, but Example \ref{exampleCleq1} shows that $\E[A_i]g_i(0)$ becomes small (or even zero) very rapidly as $\C_{A_i}^2$ gets smaller.

Summarising, the LT limit of a $GI/G/1$ queue (ignoring $\O(\rho_i^2)$ terms and higher) is:
 \begin{equation}
\E[W_{i,GI/G/1}^{\textit{LT}}]=\rho_i\, \E[A_i]g_i(0) \E[B_i^\textit{res}].\label{EWgg1}
 \end{equation}
For Poisson arrivals $(\E[A_i]g_i(0)=1)$, it is known that $\E[W_{i,M/G/1}]=\frac{\rho_i}{1-\rho_i}\E[B_i^\textit{res}] = \rho_i \E[B_i^\textit{res}]+\O(\rho_i^2)$, which is consistent with our approximation.

The second step in determining the LT limit of the mean waiting time of a type $i$ customer in a polling system, is finding the LT limits of $\E[I_i^{\textit{res}}]$, the mean residual \emph{intervisit time} of $Q_i$, and (for gated service only) $\frac{\E[V_iI_i]}{\E[I_i]}$, the mean visit time $V_i$ given that it is observed from the following intervisit time $I_i$. In this LT analysis we need to focus on first order terms only. Noting the fact that $I_i = S_i + V_{i+1}+S_{i+1}+\dots+V_{i+N-1}+S_{i+N-1}$, we condition on the moment at which $I_i$ is observed. We distinguish between two cases. The moment of observation either takes place during a visit time, or during a switch-over time:
\begin{align}
\E[I_i^{\textit{LT,res}}] &= \sum_{j=1}^{N-1}\frac{\E[V_{i+j}]}{\E[I_i]}\E[I_i^{\textit{LT,res}}|\text{observed during }V_{i+j}] \nonumber\\
&+ \sum_{j=0}^{N-1}\frac{\E[S_{i+j}]}{\E[I_i]}\E[I_i^{\textit{LT,res}}|\text{observed during }S_{i+j}].
\label{EIiresLT}
\end{align}

\paragraph{Observation during visit time.}
The probability that a random observation epoch takes place during a visit time, say $V_j$, is $\frac{\E[V_j]}{\E[I_i]}$, for any $j\neq i$. However, we are only interested in order $\rho$ terms, so this probability simplifies to
\[\frac{\E[V_j]}{\E[I_i]} = \frac{\rho_j \E[C]}{(1-\rho_i)\E[C]} = \rho_j + \O(\rho^2).\]
The fact that this probability is $\O(\rho)$, implies that all further $\O(\rho)$ terms can be ignored in $\E[I_i^{\textit{LT,res}}|\text{observed during }V_j]$, because in LT we focus on first order terms only.

The length of the residual intervisit time is the length of the residual visit period of type $j$ customers, $V_j^\textit{res}$, plus all switch-over times $S_j +\dots+S_{i-1}$, plus all visit times $V_{j+1} +\dots+V_{i-1}$.
The first term simplifies to $\E[V_j^\textit{res}] = \E[B_j^\textit{res}]+\O(\rho)$. The terms $\E[V_k|\text{observed from }V_j], k=j+1,\dots,i-1$, in light traffic, are all $\O(\rho)$. Summarising, the mean residual intervisit period when observed during $V_j$ is simply a mean residual service time $\E[B_j^\textit{res}]$, plus all mean switch-over times $\E[S_j +\dots+S_{i-1}]$, plus $\O(\rho)$ terms:
\begin{equation}
\E[I_i^{\textit{LT,res}}|\text{observed during }V_j] = \E[B_j^\textit{res}]+\sum_{k=j}^{i-1}\E[S_k] + \O(\rho).
\label{EIiduringVj}
\end{equation}
The intuition behind this equation is that, in light traffic,  the probability of having another service during the residual cycle is negligible, i.e., $\O(\rho)$. Hence, the length of the residual intervisit time is solely determined by the residual service time and the remaining switch-over times in the
cycle.


\paragraph{Observation during switch-over time.}
We continue by determining the mean residual intervisit period, conditioned on a random observation epoch during a switch-over time, say $S_j$, $j=1,\dots,N$. The probability that such an epoch takes place during $S_j$, is
\[\frac{\E[S_j]}{\E[I_i]} = \frac{\E[S_j]}{(1-\rho_i)\E[C]} = \frac{\E[S_j]}{\E[S]}\frac{1-\rho}{1-\rho_i} = \frac{\E[S_j]}{\E[S]}(1-\rho+\rho_i)+\O(\rho^2).\]
It becomes apparent from this expression that things get slightly more complicated now, because order $\rho$ terms in the conditional residual intervisit time may no longer be neglected. The residual intervisit time now consists of the residual switch-over time $S_j^\textit{res}$, plus the switch-over times $S_{j+1} +\dots+S_{i-1}$, plus all visit periods $V_{j+1} +\dots+V_{i-1}$. The length of a visit period $V_k$, for $k>j$, is the sum of the busy periods of all type $k$ customers that have arrived during $S_{k-N}, \dots, S_{j-1}$, $S_j^\textit{past}$, $S_j^\textit{res}$, and $S_{j+1}, \dots, S_{k-1}$. By $S_j^\textit{past}$ we denote the elapsed switch-over time during which the intervisit period is observed, which has the same distribution as the residual switch-over time $S_j^\textit{res}$. Compared to an observation during a visit time, it is more difficult to determine the conditional mean length of a busy period $\E[V_k|\text{observed during }S_j]$ under LT. We use a heuristic approach, which is exact if the arrival process of type $k$ customers is Poisson, and approximate it by:
\[\E[V_k|\text{observed during }S_j] \approx \rho_k\left(\sum_{l\neq j}\E[S_l] + \E[S_j^\textit{past}] + \E[S_j^\textit{res}]\right) + \O(\rho^2), \qquad k=j+1,\dots,i-1.\]
If $A_k$ is exponentially distributed, the above expression is exact.
Nevertheless, numerical experiments have shown that this approximative assumption has no or at least negligible impact on the accuracy of the approximated mean waiting times. Summarising:
\begin{align}
\E[I_i^{\textit{LT,res}}|\text{observed during }S_j]&\approx\sum_{k=i}^{j-1}\E[S_k]\big(\sum_{l=j+1}^{i+N-1}\rho_l\big)+\E\left(S_j^{\textit{past}}\right)\big(\sum_{k=j+1}^{i+N-1}\rho_k\big) \nonumber\\
&+ \E\left(S_j^{\textit{res}}\right)\big(1+\sum_{k=j+1}^{i+N-1}\rho_k\big)
+\sum_{k=j+1}^{i+N-1}\E[S_k]\big(1+\sum_{l=j+1}^{i+N-1}\rho_l\big)
+\O(\rho^2).\label{EIiduringSj}
\end{align}

The expression for $I_i^\textit{res}$ under light traffic conditions now follows from substituting \eqref{EIiduringVj} and \eqref{EIiduringSj} in \eqref{EIiresLT}. The result can be rewritten to:
\begin{align}
\E[I_i^{\textit{LT,res}}] 
%
\approx& \sum_{j=i+1}^{i+N-1}\rho_j\E[B_j^{\textit{res}}] + \sum_{j=i+1}^{i+N-1}\rho_j\sum_{k=j}^{i+N-1}\E[S_k]\nonumber\\
&+\sum_{j=i}^{i+N-1}\frac{1}{2\E[S]}\left[\E\left(S_j^{2}\right)\big(1-\rho+\rho_i+2\sum_{k=j+1}^{i+N-1}\rho_k\big)
\right]\nonumber\\
&+\frac{1}{\E[S]}\sum_{j=i}^{i+N-1}\left[\sum_{k=i}^{j-1}\E[S_j]\E[S_k]\big(\sum_{l=j+1}^{i+N-1}\rho_l\big) + \sum_{k=j+1}^{i+N-1}\E[S_j]\E[S_k]\big(1-\rho+\rho_i+\sum_{l=j+1}^{i+N-1}\rho_l\big)
\right]\nonumber\\
&+\O(\rho^2)\nonumber\\
=& \sum_{j=i+1}^{i+N-1}\rho_j\E[B_j^{\textit{res}}] + \sum_{j=i+1}^{i+N-1}\rho_j\sum_{k=j}^{i+N-1}\E[S_k]\nonumber\\
&+(1-\rho+\rho_i)\E[S^{\textit{res}}] +\frac{1}{\E[S]}\sum_{j=i}^{i+N-1}\sum_{k=i}^{i+N-1}\E[S_jS_k]\big(\sum_{l=j+1}^{i+N-1}\rho_l\big)+\O(\rho^2)\nonumber\\
=&\,(1+\rho)\E[S^\textit{res}]+\rho\E[B^{\textit{res}}] - \rho_i\E[B_i^\textit{res}]+\rho_i\big(\E[S^{\textit{res}}]-\E[S]\big) -\frac{1}{\E[S]}\sum_{j=0}^{N-1}\sum_{k=0}^j\rho_{i+k}\V[S_{i+j}]\nonumber\\
&+\O(\rho^2),\label{EIi}
\end{align}
for $i=1,\dots,N$. The last step in \eqref{EIi} follows after some straightforward (but tedious) rewriting.

The Fuhrmann-Cooper decomposition of the mean waiting time for customers in a polling system with \emph{gated} service \eqref{decompositiongated}, also requires the computation of $\frac{\E[V_iI_i]}{\E[I_i]}$ under LT conditions. Here, again, we have to resort to using a heuristic and use $\frac{\E[V_iI_i]}{\E[I_i]} = \rho_i \E[S] + \O(\rho^2)$, because this value is \emph{exact} in the case of Poisson arrivals. 
This term is the mean length of the visit time $V_i$ given that it is observed during the \emph{following} intervisit time $I_i$. The term appears because, contrary to exhaustive service, type $i$ customers arriving during $V_i$ are not served until the next cycle. However, it is easier to consider ${\E[V_iI_i]}/{\E[V_i]}$ instead, and to use the relation
\[\frac{\E[V_iI_i]}{\E[I_i]} = \frac{\E[V_iI_i]}{\E[V_i]}\times\frac{\E[V_i]}{\E[I_i]}.\]
The term ${\E[V_iI_i]}/{\E[I_i]}$ is the mean length of the intervisit time $I_i$ following $V_i$, given that it is observed during this visit time $V_i$.
Firstly, we note that
\[\frac{\E[V_i]}{\E[I_i]} = \frac{\rho_i}{1-\rho_i} = \rho_i+\O(\rho^2).\]
This implies that we can ignore all $\O(\rho)$ terms in ${\E[V_iI_i]}/{\E[V_i]}$, which means that only the switch-over times play a role,
\[\frac{\E[V_iI_i]}{\E[V_i]} = \E[S]+\O(\rho).\]
Concluding, we have
\begin{equation}
\frac{\E[V_iI_i]}{\E[I_i]} = \rho_i \E[S] + \O(\rho^2),\label{EviIi}
\end{equation}
in the case of Poisson arrivals. If the arrival process is not Poisson, this is not exact, but we use it as an approximation.

Having made all required preparations, we are ready to formulate the main result of the present subsection.
Under light traffic, an approximation for the mean waiting time of a type $i$ customer in a polling model with general arrivals and respectively exhaustive and gated service in $Q_i$, is:
\begin{align}
\E[W_i^{\textit{LT,exh}}] \approx\,& \E[S^{\textit{res}}] + \rho_i (\E[\hat{A}_i]\hat{g}_i(0) -1)\E[B_i^\textit{res}] + \rho\E[B^{\textit{res}}] + (\rho-\rho_i)\left(\E[S]-\E[S^{\textit{res}}]\right)\nonumber\\
&+\frac{1}{\E[S]}\sum_{k=i+1}^{i+N-1}\rho_{k}\sum_{j=i}^{k-1}\V[S_j]+\O(\rho^2),\qquad i = 1,\dots,N,\label{EWltexhaustive}\\
\E[W_i^{\textit{LT,gated}}] \approx\,& \E[W_i^{\textit{LT,exh}}] + \rho_i \E[S],\label{EWltgated}
\end{align}
where $\hat{g}_i(t)$ is the density of the interarrival times of type $i$ customers at $\rho=1$. Equation \eqref{EWltexhaustive} follows from substitution of \eqref{EWgg1} and \eqref{EIi} in
\begin{equation}
\E[W_i] \approx \E[W_{i,GI/G/1}]+\E[I_i^{\textit{res}}],\qquad{i=1,\dots,N}.\label{fuhrmanncooperdecomposition}
\end{equation}

For Poisson arrivals, \eqref{EWltexhaustive} and \eqref{EWltgated} are exact. The LT limit for polling systems with Bernoulli service (and Poisson arrivals) has been experimentally found in \cite{blancvdmei95} and, indeed, it can be shown that their result for exhaustive service, which is a special case of Bernoulli service, agrees with our result after substituting $\E[\hat{A}_i]\hat{g}_i(0)=1$ in \eqref{EWltexhaustive}.

\subsection{Heavy traffic}\label{heavytraffic}

Heavy traffic limits in polling systems have been studied by Coffman et al. \cite{coffman95,coffman98}, and by Olsen and Van der Mei \cite{olsenvdmei03,olsenvdmei05}. In these papers, the HT limits of the waiting time distributions are found under the assumption of Poisson arrivals. For general renewal arrivals, a proof is given for the special case $N=2$ (cf. \cite{coffman95,coffman98}), and a strong conjecture for larger values of $N$ (cf. \cite{olsenvdmei05}).
In \cite{vdmeiwinands08}, the following result for the \emph{mean} waiting time is proven rigorously for polling systems with renewal arrivals:
\begin{equation}
\E[W_{i}^{\textit{HT}}] = \frac{\omega_i}{1-\rho} + {o}((1-\rho)^{-1}), \qquad \rho \uparrow 1.\label{EWht}
\end{equation}
Obviously, in HT, all queues become unstable and, thus, $\E[W_i]$ tends to infinity for all $i$. The rate at which $\E[W_{i}]$ tends to infinity as $\rho \uparrow 1$ is indicated by $\omega_i$, which is referred to as the \emph{mean asymptotic scaled delay} at queue $i$, and depends on the service discipline. For exhaustive service,
\[
\omega_i =  \frac{1-\hat{\rho}_i}{2} \left( \frac{\sigma^2}{\sum_{j=1}^N \hat{\rho}_j (1-\hat{\rho}_j)} +  \E[S] \right),\qquad i=1,\dots,N,
\]
with
\begin{equation*}
\sigma^2:= \sum_{i=1}^N \hat{\lambda}_i \left( \V[B_i] + \hat{\rho}_i^2 \V[\hat{A}_i] \right). \label{sigma}
\end{equation*}
Here, the limits are taken such that the arrival rates are increased, while keeping the service-time distributions fixed, and keeping the distributions of the interarrival times $A_i~(i=1,\ldots,N)$ fixed up to a common scaling constant $\rho$. Notice that in the case of Poisson arrivals we have $\sigma^2=\E[B^2]/\E[B]$.

For gated service, we have
\[
\omega_i =  \frac{1+\hat{\rho}_i}{2} \left( \frac{\sigma^2}{\sum_{j=1}^N \hat{\rho}_j (1+\hat{\rho}_j)} +  \E[S] \right).
\]

\subsection{Interpolation}\label{inter}

Now that we have the expressions for the mean delay in both LT and HT, we can determine the constants $K_{0,i}, K_{1,i}$, and $K_{2,i}$ in approximation formula \eqref{ewapprox}. We simply impose the requirements that approximation \eqref{ewapprox} results in the same mean waiting time for $\rho=0$ as the LT limit, and for $\rho\uparrow1$ as the HT limit. Since \eqref{EWltexhaustive} (and \eqref{EWltgated} for gated service) has been determined up to the first order of $\rho$ terms, we also add the requirement that the derivative with respect to $\rho$, taken at $\rho=0$, of our approximation is equal to the derivative of the LT limit.
A more formal definition of these requirements is presented below:
\begin{align*}
\E[W_{i,\textit{app}}]\big|_{\rho=0} &= \E[W_{i}]\big|_{\rho=0},\\
\frac{\dd}{\dd\rho}\E[W_{i,\textit{app}}]\big|_{\rho=0} &= \frac{\dd}{\dd\rho}\E[W_{i}]\big|_{\rho=0},\\
(1-\rho)\E[W_{i,\textit{app}}]\big|_{\rho=1} &= (1-\rho)\E[W_{i}]\big|_{\rho=1}.
\end{align*}
This leads to \eqref{ewapprox} as approximation for $\E[W_i]$ in a polling system with general arrivals.
Constants $K_{0,i}$, $K_{1,i}$, and $K_{2,i}$ are defined in \eqref{c0}--\eqref{c2} for systems with exhaustive service, or \eqref{c0gated}--\eqref{c2gated} for gated service.

\subsection{Matching properties}\label{specialcases}

A desirable property of an approximation is that it matches known exact results.
In the present section we discuss several cases where the interpolation approximation yields \emph{exact} results. Most cases require Poisson arrivals, but it is shown that also in
two limiting cases where exact results are available for general arrivals, i.e., heavy traffic and large switch-over times, the approximation is exact.
It further turns out that, in case of Poisson arrivals, the approximated mean waiting times satisfy the \emph{pseudo-conservation law}, implying that the weighted sum
$\sum_{i=1}^N \rho_i \E[W_{i,\textit{app}}]$ is exact for each load $0 \le \rho <1$. In fact, this appears to be true for any interpolation approximation for the mean waiting times of the form \eqref{inter1} matching the HT limit and the LT limits of order $0$ up to order $k$, provided $k > 0$.
These properties of the interpolation approximation indicate that it is the ``natural'' approximation, given the HT and LT limits.

\paragraph{Light and heavy traffic.}

The light traffic limit of $\E[W_i]$, given by \eqref{EWltexhaustive} for exhaustive service and by \eqref{EWltgated} for gated service, is exact for Poisson arrivals. The heavy traffic limit \eqref{EWht} of $\E[W_i]$ is even exact for renewal arrivals. An appropriate choice of constants $K_{0,i}, K_{1,i}$ and $K_{2,i}$ can reduce \eqref{ewapprox} to either \eqref{EWltexhaustive}, \eqref{EWltgated}, or \eqref{EWht}. Since the LT and HT limits have been used in the set of equations that determine the coefficients of the approximation, it goes without saying that $\E[W_{i,\textit{app}}]$ is equal to \eqref{EWltexhaustive} (or \eqref{EWltgated} for gated service) and \eqref{EWht}, for $\rho \downarrow 0$ and $\rho \uparrow 1$ respectively. This implies that the LT limit of our approximation is exact for Poisson arrivals, and the HT limit is exact for general arrivals.

\paragraph{Symmetric system.}

If $\hat\rho_i=\frac1N$ for all $i=1,\dots,N$, all $B_i$ have the same distribution, and the variances $\V[S_i]$ of all switch-over times are equal, then our approximation is exact if all interarrival distributions are exponential. For exhaustive service, we obtain
\begin{align*}
K_{1,i} = \,&\E[B^{\textit{res}}] + \frac{N-1}{N}\E[S]-\left(2-\frac1N\right)\E[S^{\textit{res}}]+\frac{1}{\E[S]}\sum_{k=i+1}^{i+N-1}\hat\rho_{k}\sum_{j=i}^{k-1}\V[S_j]\\
= \,&\E[B^{\textit{res}}] + \frac{N-1}{N}\E[S]-\left(2-\frac1N\right)\E[S^{\textit{res}}]+\frac{N-1}{N}\frac{\V[S]}{2\E[S]}\\
 = \,&\E[B^{\textit{res}}] + \left(1-\frac1N\right)\frac{\E[S]}{2} -\E[S^\textit{res}],
 \end{align*}
and $K_{2,i} = 0$. This implies that $\E[W_{i,\textit{app}}] = \E[W_{i,\textit{symm}}]$, since
\[
\E[W_{i,\textit{symm}}] =
\frac{\rho}{1-\rho} \E[B^\textit{res}] + \E[S^\textit{res}] +\frac{\rho(1-\frac{1}{N})}{1-\rho}\frac{\E[S]}{2}.
\]
Note that $\E[W_{i,\textit{symm}}]$ is of the form (\ref{inter1}) with $q(\rho)$ being a linear polynomial. In fact, this immediately implies that the interpolation is exact, given the (exact) LT and HT limits. Further, we do \emph{not} require that the mean switch-over times $\E[S_i]$ are equal. One can verify that the same holds for gated service.

\paragraph{Single queue (vacation model).}

An immediate consequence of the fact that our approximation is exact in symmetric polling systems with Poisson arrivals, is that it also gives exact results for the mean waiting time of customers in a single-queue polling system with Poisson arrivals. A polling system consisting of only one queue, but with a switch-over time between successive visits to this queue, is generally referred to as a queueing system with multiple server vacations.

\paragraph{Large switch-over times.}

For $S$ deterministic, $S \rightarrow \infty$, and, again, under the assumption of Poisson arrivals, it is proven in \cite{winandslargesetups07,winandslargesetupsbranching09} that $\frac{\E[W_{i}]}{S} \rightarrow \frac{1-\rho_i}{2(1-\rho)}$ for exhaustive service. It can easily be verified that our approximation has the same limiting behaviour:
\begin{equation}
\lim_{S\rightarrow\infty} \frac{\E[W_{i,\textit{app}}]}{S} = \frac{1-\rho_i}{2(1-\rho)}.\label{Sinfty}
\end{equation}
For gated service, $\frac{\E[W_{i,\textit{app}}]}{S} \rightarrow \frac{1+\rho_i}{2(1-\rho)}$, which is also the exact limit (see, e.g., \cite{winandsPhD}).

Olsen \cite{olsen2001} studies the effects of large switch-over times in polling systems operating under HT conditions. She discovers that, under these conditions, the limiting behaviour \eqref{Sinfty} is also exhibited in polling systems with \emph{general renewal arrivals}.

\paragraph{Miscellaneous other exact results.}

The approximation is also exact in several other cases, all with Poisson arrivals, when the parameter values are carefully chosen. The relations between the input parameters that yield exact approximation results become very complicated, especially in polling systems with more than two queues. We only mention one interesting example here: our approximation gives exact results for a two-queue polling system with exhaustive service and
\begin{equation}
\E[B_1] = \E[B_2], \E[S_1] = \E[S_2], \C_{A_1}^2 = \C_{A_2}^2, \C_{B_1}^2 = \C_{B_2}^2, \C_{S_1}^2 = \C_{S_2}^2,\label{miscexact1}
\end{equation}
if the following constraint is satisfied:
\begin{equation}
\rho = \frac{1+I_{A_i}^2}{2I_{A_i}}-\frac{\C_{S_i}^2}{1+\C_{B_i}^2}\cdot\frac{\E[S_i]}{\E[B_i]},\label{miscexact2}
\end{equation}
where $I_{A_i} = \frac{\hat\rho_1}{\hat\rho_2}$ is the ratio of the loads of the two queues.
Obviously, if $I_{A_i} = 1$, the system is symmetric and our approximation gives exact results regardless of the other parameter settings.

\paragraph{Pseudo-conservation law.}

A well-known result in polling literature, is the so-called \emph{pseudo-conservation law}, derived by Boxma and Groenendijk \cite{boxmagroenendijk87} using the concept of work decomposition. This law gives the following exact expression for the weighted sum of the mean waiting times in a polling system with Poisson arrivals:
\begin{equation}
\sum_{i=1}^N \rho_i \E[W_i] = \frac{\rho^2}{1-\rho}\E[B^\textit{res}]+\rho \E[S^\textit{res}]+\frac{\E[S]}{2}\frac{\rho^2-\sum_{j=1}^N\rho_j^2}{1-\rho}+\sum_{j=1}^N\E[Z_{jj}],\label{pcl}
\end{equation}
where $\E[Z_{jj}]$ denotes the mean amount of work left behind in $Q_j$ at the completion of a visit of the server to $Q_j$. It is shown in \cite{boxmagroenendijk87} that $\E[Z_{jj}]$ is the only term that depends on the service discipline. For exhaustive service $\E[Z_{jj}] = 0$, for gated service $\E[Z_{jj}] = \rho_j^2\frac{\E[S]}{1-\rho}$.
It can be shown that the interpolation approximation \eqref{ewapprox} satisfies the pseudo-conservation law in the case of Poisson arrivals: if $\E[\hat{A}_i]\hat{g}_i(0) = 1$ for $i=1,\dots,N$, then $\sum_{i=1}^N \rho_i \E[W_{i,\textit{app}}]$ equals the right-hand side of \eqref{pcl}.
In fact, in case of Poisson arrivals, {\em any interpolation approximation} for the mean waiting times of the form \eqref{inter1}, matching the HT limit and the LT limits of order $0$ up to order $k$, satisfies \eqref{pcl}, provided $k > 0$. To establish this result, we first rewrite \eqref{pcl} in the form $Q(\rho)=0$, by moving the terms at right-hand side of \eqref{pcl} to the left. Now let $\tilde{Q}(\rho)$ denote the version of $Q(\rho)$ where the mean waiting times have been replaced by their interpolation approximation. Clearly, the interpolation approximation and the right-hand side of \eqref{pcl} are {\em both} of the form \eqref{inter1}, and thus $\tilde{Q}(\rho)$ is also of the form
\[
\tilde{Q}(\rho) = \frac{\tilde{q}(\rho)}{1-\rho},
\]
where $\tilde{q}(\rho)$ is a polynomial of degree $k+2$. Since the interpolation approximation matches the LT limits of order $0$ up to order $k$, it follows that $\tilde{Q}^{(n)} (0) = Q^{(n)} (0) = 0$, and hence, $\tilde{q}^{(n)} (0) = 0$ for all $n = 0, 1, \ldots, k+1$. Further, by the HT limit,
$\tilde{q}(1) = \lim_{\rho \rightarrow 1} (1-\rho) \tilde{Q} (\rho) = \lim_{\rho \rightarrow 1} (1-\rho) Q (\rho) = 0$. This implies that $\tilde{q}(\rho) = 0$, and thus $\tilde{Q} (\rho) = 0$ for all $\rho$.

\section{Numerical study}\label{numericalexample}

\subsection{Initial glance at the approximation}\label{firstexample}

\begin{figure}[h!t]
\begin{center}
\includegraphics[width=0.8\linewidth]{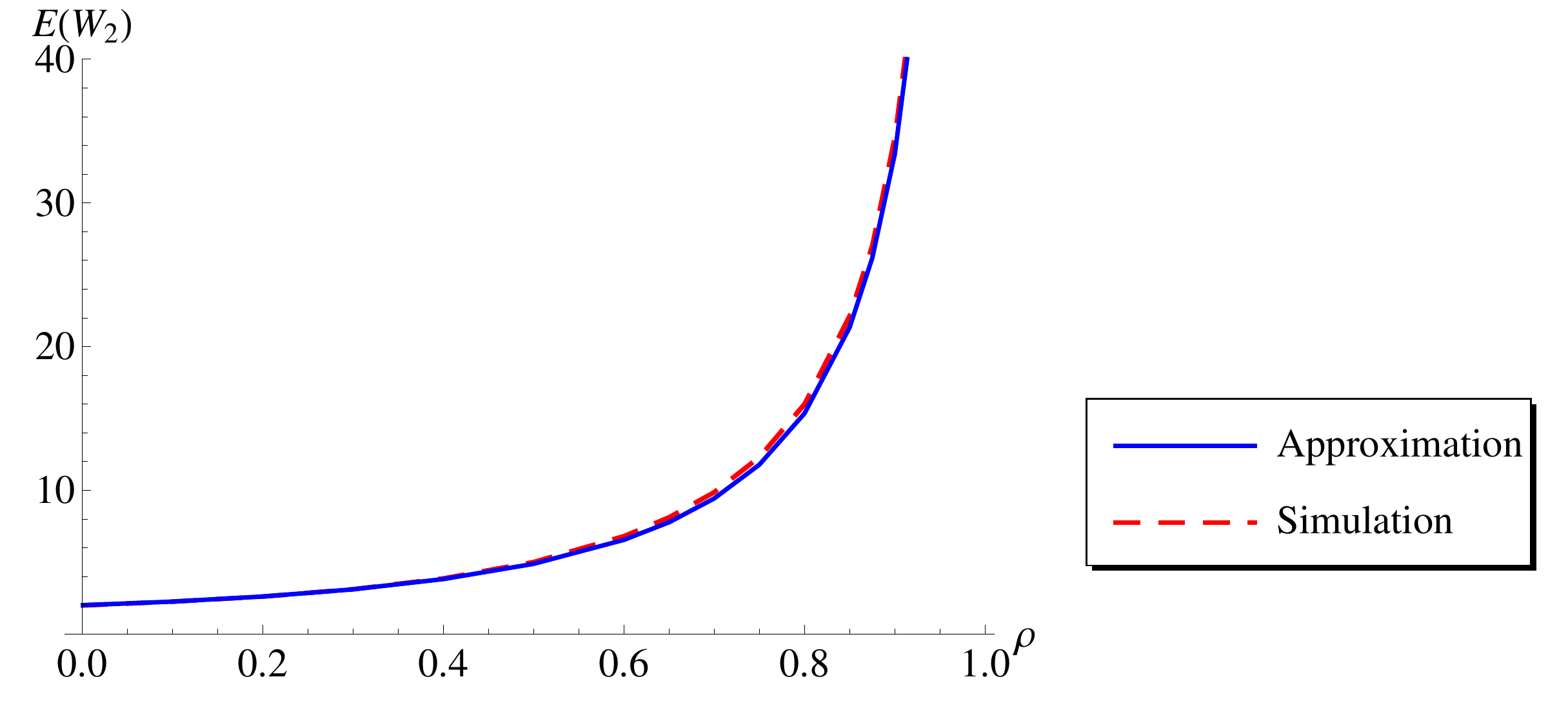}
\end{center}
\caption{Approximated and simulated mean waiting time $\E[W_2]$ of $Q_2$ of the example in subsection \ref{firstexample}.}
\label{figExample1}
\end{figure}

Before we study the accuracy of the approximation to a large test bed of polling systems, we just pick a rather arbitrary, simple system to compare the approximation with exact results in order to get some initial insights. Consider a three-queue polling system with loads of $Q_1$, $Q_2$, and $Q_3$ divided as follows: $\hat\rho_1 = 0.1, \hat\rho_2 = 0.3$, and $\hat\rho_3 = 0.6$. All service times and switch-over times are exponentially distributed, with mean 1. The interarrival times have SCV $\C_{A_i}^2 = 3$ for $i=1,2,3$. In Figure \ref{figExample1} we plot the approximated mean waiting time of $Q_2$, $\E[W_{2,\textit{app}}]$, versus the load of the system $\rho$. Since this system cannot be analysed analytically, we compare the approximated values with simulated values.
Both in the approximation and in the simulation we fit a $H_2$ distribution as described in Example \ref{exampleCgeq1}.

The errors are largest for $Q_2$, which is the reason why we chose this queue in particular in Figure \ref{figExample1}. The most important information that this figure reveals, is that even though the accuracy of the approximation is worst for this queue (a relative error of $-4.47\%$ for $\rho=0.7$), the shape of the approximation function is very close to the shape of the exact function, which makes it very suitable for optimisation purposes. The maximum relative errors of $Q_1$ and $Q_3$ are $3.10\%$ and $2.90\%$ respectively.

In order to get more insight in the numerical accuracy of the approximation for a wide variety of different parameter settings, we create a large test bed in the next subsection and compare the approximation with exact or simulated results. It turns out that the maximum relative errors for most of the polling systems are smaller than the one selected in the above example.

\subsection{Accuracy of the approximation}

In the present section we study the accuracy of our approximation. We compare the approximated mean waiting times of customers in various polling systems to the exact values. The complete test bed of polling systems that are analysed, contains 2304 different combinations of parameter values, all listed in Table \ref{testbed}. We show detailed results for exhaustive service first, and discuss polling systems with gated service at the end of this section.
\begin{table}[h!t]
\begin{center}
\begin{tabular}{|l|c|c|}
\hline
Parameter & Notation & Values\\
\hline
Number of queues & $N$ & $2, 3, 4, 5$ \\
Load & $\rho$ & $0.1, 0.3, 0.5, 0.7, 0.9, 0.99$ \\
SCV interarrival times & $\C^2_{A_i}$ & $0.25, 1, 2$ \\
SCV service times & $\C^2_{B_i}$ & $0.25, 1$ \\
SCV switch-over times & $\C^2_{S_i}$ & $0.25, 1$ \\
Imbalance interarrival times & $I_{A_i}$ & $1, 5$ \\
Imbalance service times & $I_{B_i}$ & $1, 5$ \\
Ratio service and switch-over times & $I_{S_i/B_i}$ & $1, 5$ \\
\hline
\end{tabular}
\end{center}
\caption{Test bed used to compare the approximation to exact results.}
\label{testbed}
\end{table}
We have varied the load between $0.1$ and $0.9$ with steps of $0.2$, and included $\rho=0.99$ to analyse the limiting behaviour of our approximation when the load tends to $1$.
The SCV of the interarrival times, $\C^2_{A_i}$, is varied between $0.25$ and $2$. In case of non-Poisson arrivals, i.e. $\C^2_{A_i} \neq 1$, the exact values have been established through extensive simulation because they cannot be obtained in an analytic way. In these simulations we fit a phase-type distribution to the first two moments of the interarrival times, as described in Examples \ref{exampleCgeq1} and \ref{exampleCleq1}. For service times and switch-over times, only SCVs of $0.25$ and $1$ are considered. SCVs greater than $1$ are less common in practice and are discussed separately from the test bed later in this section. The imbalance in interarrival times and service times, $I_{A_i}$ and $I_{B_i}$, is the ratio between the largest and the smallest mean interarrival/service time.
The interarrival times are determined in such a way, that the overall mean is always $1$, $\lambda_1$ is the largest and $\lambda_N$ the smallest, and the steps between the $\lambda_i$ are linear. E.g., for $N=5$ and $I_{A_i} = 5$ we get $\lambda_i = 2-i/3, i=1,\dots,5$. The mean service times $\E[B_i]$ increase linearly in $i=1,\dots,N$, with $\E[B_N] = I_{B_i}\E[B_1]$ (so $\E[B_1]$ is the smallest mean service time). They follow from the relation $\sum_{i=1}^N \lambda_i\E[B_i]=\rho$. E.g., for $N=5$, and $I_{A_i} = I_{B_i} = 5$ we get $\E[B_i]/ \rho = 3i/35$. The last parameter that is varied in the test bed, is the ratio between the mean switch-over times and the mean service times, $I_{S_i/B_i} = \frac{\E[S_i]}{\E[B_i]}$. The total number of systems analysed is $4\times6\times3\times2^5=2304$. A system consisting of $N$ queues results $N$ mean waiting times, $\E[W_1], \dots, \E[W_N]$, so in total these $2304$ systems yield $8064$ mean waiting times. The absolute relative errors, defined as $|o-e|/e$, where $o$ stands for observed (approximated) value, and $e$ stands for expected (exact) value, are computed for all these $8064$ queues. Table \ref{numericalresults1} shows these relative errors (times $100\%$) categorised in bins of $5\%$. In this table, and in all other tables, results for systems with a different number of queues are displayed in separate rows. The reader should keep in mind that the statistics in each row are based on $\frac14\times2304\times N$ absolute relative errors, where $N$ is the number of queues used in the specified row.
Table \ref{numericalresults1} shows that, e.g., $98.84\%$ of the approximated mean waiting times in polling systems consisting of $3$ queues deviate less than $5\%$ from their true values.
From Table \ref{numericalresults1} it can be concluded that the approximation accuracy increases with the number of queues in a polling system.
More specifically, for systems with more than 2 queues, no approximation errors are greater than 10\%, and the vast majority is less than 5\%. The mean relative errors for $N=2,\dots,5$ are respectively $2.18\%, 0.93\%, 0.70\%$ and $0.57\%$. 
It is also noteworthy, that $193$ out of the $2304$ systems yield exact results. All of these $193$ systems have Poisson input, and all of them -- except for one -- are symmetric. The only asymmetric case for which our approximation yields an exact result, happens to satisfy constraints \eqref{miscexact1} and \eqref{miscexact2}.

In Table \ref{errorpervar} the mean relative error percentages are shown for a combination of input parameter settings. The number of queues is always varied per row, while per column another input parameter is varied. 
This way we can find in more detail which (combinations of) parameter settings result in large approximation errors. In Table \ref{errorpervar}(a) the load $\rho$ is varied, and it can be seen that for a load of $\rho=0.7$ the approximation is least accurate. E.g., the mean relative error of all approximated waiting times in polling systems consisting of 3 queues with a load of $\rho=0.7$ is $1.69\%$. Table \ref{errorpervar}(b) shows the impact of the SCV of the interarrival times on the accuracy. Especially for systems with more than 2 queues the accuracy is very satisfactory, in particular for the case $\C^2_{A_i} = 1$. In Table \ref{errorpervar}(c) the impact of imbalance in a polling system on the accuracy is depicted, and, as could be expected, it can be concluded that a high imbalance in either service or interarrival times has a considerable, negative, impact on the approximation accuracy. Polling systems with more than 2 queues are much less bothered by this imbalance than polling systems with only 2 queues. 

\begin{table}[h!t]
\[
\begin{array}{|c|cccc|}
\hline
N  & 0-5\% & 5-10\% & 10-15\% & 15-20\% \\
\hline
2 &   86.46 & 10.24 & 2.78 & 0.52 \\
3 &   98.84 & 1.16 & 0.00 & 0.00 \\
4 &   99.78 & 0.22 & 0.00 & 0.00 \\
5 &   99.93 & 0.07 & 0.00 & 0.00 \\
\hline
\end{array}
\]
\caption{Errors of the approximation applied to the 2304 test cases with exhaustive service, as described in Section \ref{numericalexample}, categorised in bins of $5\%$.}
\label{numericalresults1}
\end{table}

\begin{table}[h!t]
\begin{center}
\parbox{0.45\textwidth}{
\begin{center}
\[
\begin{array}{|c|cccccc|}
\hline
N & \multicolumn{6}{|c|}{\textrm{Load ($\rho$)}} \\
\hline
  & 0.10 & 0.30 & 0.50 & 0.70 & 0.90 & 0.99\\
\hline
2 &  0.31 & 1.81 & 3.41 & 4.17 & 2.70 & 0.67 \\
3 &  0.16 & 0.84 & 1.44 & 1.69 & 1.07 & 0.39 \\
4 &  0.13 & 0.68 & 1.14 & 1.28 & 0.73 & 0.25 \\
5 &  0.11 & 0.57 & 0.94 & 1.03 & 0.57 & 0.22 \\
\hline
\end{array}
\]
(a)
\end{center}
}
\hfill
\parbox{0.45\textwidth}{
\begin{center}
\[
\begin{array}{|c|ccccccccc|}
\hline
N & \multicolumn{9}{|c|}{\textrm{SCV interarrival times ($\C^2_{A_i}$)}} \\
\hline
  &&&  0.25 && 1    && 2    &&\\
\hline
2 &&&   2.27 && 1.76 && 2.50 &&\\
3 &&&   1.36 && 0.52 && 0.92 &&\\
4 &&&   1.13 && 0.29 && 0.69 &&\\
5 &&&   0.97 && 0.19 && 0.56 &&\\
\hline
\end{array}
\]
(b)
\end{center}
}
\[
\begin{array}{|c|cccc|}
\hline
N & \multicolumn{4}{|c|}{\textrm{Imbalance interarrival and service times}} \\
\hline
  & I_{A_i} = 1, I_{B_i} = 1 & I_{A_i} = 1, I_{B_i} = 5 & I_{A_i} = 5, I_{B_i} = 1 & I_{A_i} = 5, I_{B_i} = 5\\
\hline
2 &   0.69 & 2.92 & 2.80 & 2.30 \\
3 &   0.65 & 1.27 & 0.75 & 1.06 \\
4 &   0.56 & 0.89 & 0.62 & 0.73 \\
5 &   0.49 & 0.69 & 0.53 & 0.59 \\
\hline
\end{array}
\]
(c)
\end{center}
\caption{Mean relative approximation error, categorised by number of queues $(N)$ and total load of the system (a), SCV interarrival times (b), and imbalance of the interarrival and service times (c).}
\label{errorpervar}
\end{table}

\subsection{Miscellaneous other cases}

In this subsection we discuss several cases that are left out of the test bed because they might not give any new insights, or because the combination of parameter values might be rarely found in practice.

\paragraph{More queues.}
Firstly, we discuss polling systems with more than 5 queues briefly. Without listing the actual results, we mention here that the approximations become more and more accurate when letting $N$ grow larger, and still varying the other parameters in the same way as is described in Table \ref{testbed}. For $N=10$ already, all relative errors are less than $5\%$, with an average of less than $0.5\%$  and it only gets smaller as $N$ grows further.

\paragraph{More variation in service times and switch-over times.}
In the test bed we only use SCVs $0.25$ and $1$ for the service times and switch-over times, because these seem more relevant from a practical point of view. As the coefficient of variation grows larger, our approximation will become less accurate. E.g., for Poisson arrivals we took $\C^2_{B_i} \in \{2,5\}$, $\C^2_{S_i} \in \{2,5\}$, and varied the other parameters as in our test bed (see Table \ref{testbed}). This way we reproduced Table \ref{numericalresults1}. The result is shown in Table \ref{numericalresults1highSCVs} and indicates that the quality of our approximation deteriorates in these extreme cases. The mean relative errors for $N=2,\dots,5$ are respectively $3.58\%$, $1.78\%$, $1.07\%$, and $0.77\%$, which is still very good for systems with such high variation in service times and switch-over times. For non-Poisson input, no investigations were carried out because the results are expected to show the same kind of behaviour.
\begin{table}[h!t]
\[
\begin{array}{|c|cccc|}
\hline
N  & 0-5\% & 5-10\% & 10-15\% & 15-20\% \\
\hline
2 &   74.22 & 14.84 & 6.51 & 2.08 \\
3 &   89.76 & 7.29 & 2.08 & 0.69 \\
4 &   94.53 & 4.56 & 0.91 & 0.00 \\
5 &   97.71 & 2.19 & 0.10 & 0.00\\
\hline
\end{array}
\]
\caption{Errors of the approximation applied to the 768 test cases with Poisson arrival processes and high SCVs of the service times and switch-over times, categorised in bins of $5\%$.}
\label{numericalresults1highSCVs}
\end{table}

\begin{figure}[h!t]
\includegraphics[height=4.2cm]{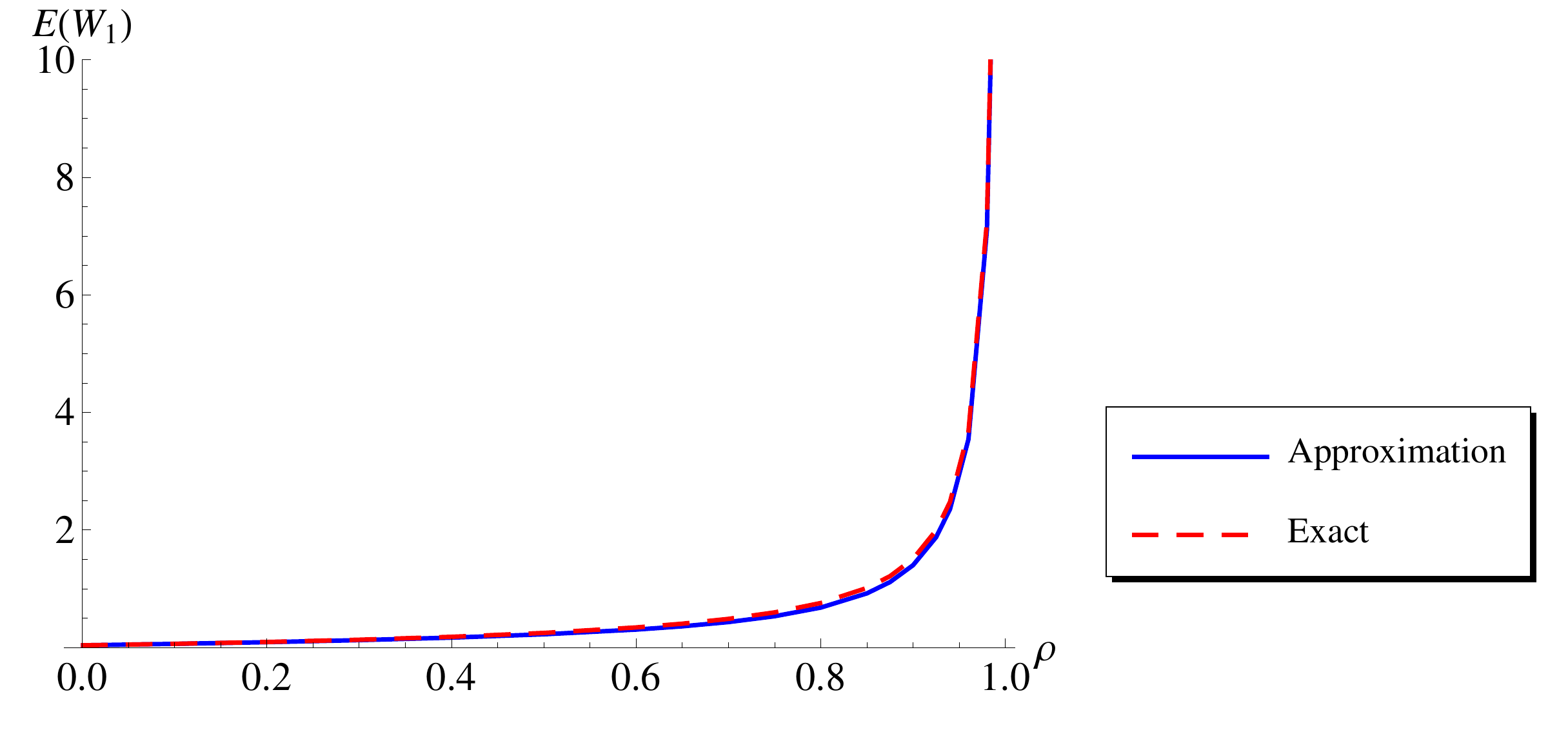}
\hfill
\includegraphics[height=4.2cm]{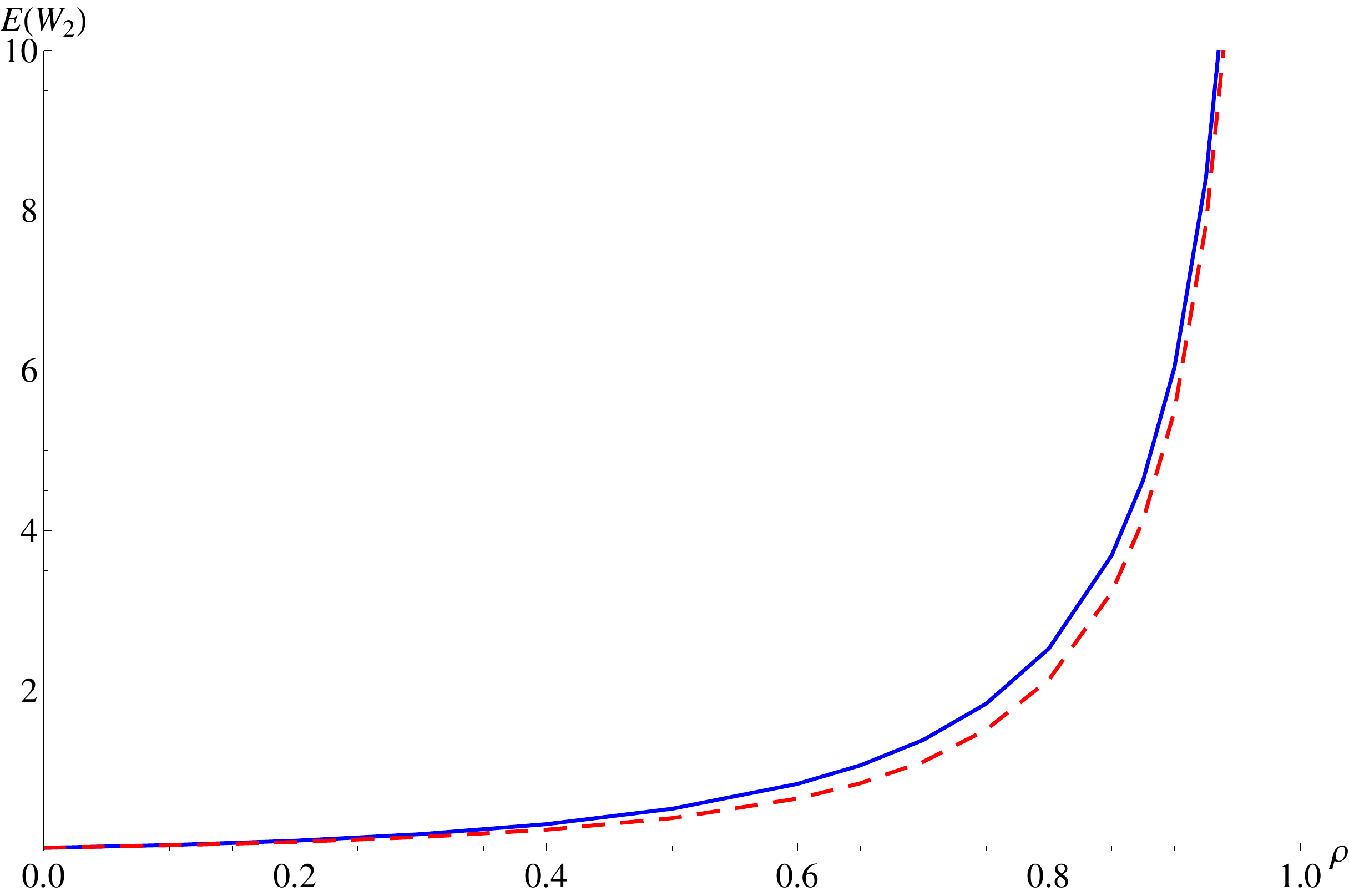}
\caption{Approximated and exact mean waiting times for a two-queue polling system with small switch-over times.}
\label{extreme3}
\end{figure}

\paragraph{Small switch-over times.}
Systems with small switch-over times, in particular smaller than the mean service times, also show a deterioration of approximation accuracy - especially in systems with 2 queues. In Figure \ref{extreme3} we show an extreme case with $N=2$, service times and switch-over times are exponentially distributed with $\E[B_i]=\frac{9}{40}$ and $\E[S_i]=\frac{9}{200}$ for $i=1,2$, which makes the mean switch-over times 5 times \emph{smaller} than the mean service times. Furthermore, the interarrival times are exponentially distributed with $\lambda_1 = 5\lambda_2$. In Figure \ref{extreme3} the mean waiting times of customers in both queues are plotted versus the load of the system. Both the approximation and the exact values are plotted. For customers in $Q_1$ the mean waiting time approximations underestimate the true values, which leads to a maximum relative error of $-11.2\%$ for $\rho=0.7$ ($\E[W_{1,\textit{app}}] = 0.43$, whereas $\E[W_1] = 0.49$). The approximated mean waiting time for customers in $Q_2$ is systematically overestimating the true value. The maximum relative error is attained at $\rho=0.5$ and is $28.8\%$ ($\E[W_{1,\textit{app}}] = 0.41$, whereas $\E[W_1] = 0.52$). Although the relative errors are high in this situation, the absolute errors are still rather small compared to the mean service time of an individual customer. This implies that the mean \emph{sojourn time} is already much better approximated. Nevertheless, this example illustrates one of the situations where our approximation gives unsatisfactory results.

\subsection{Comparison with existing approximations}

The only alternative approximations that exist for polling systems with non-exponential interarrival times, only perform well under extreme, limiting conditions. In \cite{olsenvdmei05,vdmeiwinands08} it is suggested to use the HT limit \eqref{EWht} as an approximation, but the accuracy is only found to be acceptable for $\rho > 0.8$. Another approximation for the mean waiting time in polling systems with non-exponential interarrival times uses the limit for $S\rightarrow\infty$ \cite{winandslargesetups07,winandslargesetupsbranching09}. This approximation is usable if either the total switch-over time in the system is large and the switch-over times have low variance, or if the total switch-over time in the system is large and the system is in heavy traffic. For completeness, we mention that it is also possible to construct an approximation that is purely based on the LT limit, developed in the present paper:
\begin{equation}
\E[W_{i}^{\textit{LT}}] \approx \frac{K_{0,i}+(K_{1,i} - K_{0,i})\rho}{1-\rho}.\label{EWltapprox}
\end{equation}
We do not wish to go into further details on this topic, because the accuracy of \eqref{EWltapprox} turns out to be worse for high loads. This makes our approximation, which is exact in all these limiting cases, the only one which can be applied under all circumstances. We support this statement by reproducing Table \ref{errorpervar}(a) for the three alternative approximations, based on HT limit, large switch-over times, and LT limit. A comparison of Table \ref{errorpervarHT}, which displays these results, to Table \ref{errorpervar}(a), clearly indicates the drawbacks of using approximations based on one limiting case only.
\begin{table}[h!]
\begin{center}
\[
\begin{array}{|c|cccccc|}
\hline
N & \multicolumn{6}{|c|}{\textrm{Load ($\rho$)}} \\
\hline
  & 0.10 & 0.30 & 0.50 & 0.70 & 0.90 & 0.99\\
\hline
2 &   58.97 & 51.03 & 40.34 & 27.47 & 11.34 & 1.46 \\
3 &   30.62 & 25.66 & 19.63 & 12.61 & 4.56 & 0.62 \\
4 &   22.49 & 18.61 & 14.02 & 8.80 & 3.08 & 0.42 \\
5 &   18.06 & 14.83 & 11.10 & 6.88 & 2.36 & 0.37\\
\hline
\end{array}
\]
(a)
\[
\begin{array}{|c|cccccc|}
\hline
N & \multicolumn{6}{|c|}{\textrm{Load ($\rho$)}} \\
\hline
  & 0.10 & 0.30 & 0.50 & 0.70 & 0.90 & 0.99\\
\hline
2 &  27.13 & 31.08 & 35.95 & 41.32 & 46.78 & 49.23 \\
3 &  19.28 & 21.68 & 24.51 & 27.52 & 30.33 & 31.44 \\
4 &  15.05 & 16.75 & 18.76 & 20.92 & 22.93 & 23.67 \\
5 &  12.37 & 13.67 & 15.24 & 16.96 & 18.56 & 19.17 \\
\hline
\end{array}
\]
(b)
\[
\begin{array}{|c|cccccc|}
\hline
N & \multicolumn{6}{|c|}{\textrm{Load ($\rho$)}} \\
\hline
  & 0.10 & 0.30 & 0.50 & 0.70 & 0.90 & 0.99\\
\hline
2 &    0.34 & 2.65 & 7.00 & 13.42 & 21.90 & 26.82 \\
3 &    0.24 & 1.68 & 3.90 & 6.70 & 9.62 & 10.82 \\
4 &    0.19 & 1.31 & 2.97 & 4.87 & 6.75 & 7.46 \\
5 &    0.16 & 1.09 & 2.45 & 4.02 & 5.54 & 6.12 \\
\hline
\end{array}
\]
(c)
\end{center}
\caption{Mean relative approximation error, categorised by number of queues $(N)$ and total load of the system, for three alternative approximations: based on the HT limit (a), based on large switch-over times (b), and based on the LT limit (c).}
\label{errorpervarHT}
\end{table}

For polling systems with Poisson arrivals, several alternative approximations have been developed in existing literature. The best one among them (see, e.g., \cite{boxmaworkloadsandwaitingtimes89,everitt86,groenendijk89}) uses the relation $\E[W_i] = (1\pm\rho_i)\E[C_i^\textit{res}]$, where $C_i$ is the cycle time, starting at a visit \emph{completion} to $Q_i$ when service is exhaustive, and starting at a visit \emph{beginning} for gated service. By $\pm$ we mean $-$ for exhaustive service, and $+$ for gated service. The mean residual cycle time, $\E[C_i^\textit{res}]$, is assumed to be equal for all queues, i.e. $\E[C_i^\textit{res}] \approx \E[C^\textit{res}]$, and can be found by substituting $\E[W_i] \approx (1\pm\rho_i)\E[C^\textit{res}]$ in the pseudo-conservation law \eqref{pcl}. We have used this PCL-based approximation to estimate the mean waiting times of all queues in the test bed described in Table \ref{testbed}, but taking only the $768$ cases where $\C^2_{A_i}= 1$. Table \ref{numericalresults1poisson} shows the mean relative errors for our approximation (a) and the PCL approximation (b), categorised in bins of $5\%$ as was done before in Table \ref{numericalresults1}. From these tables (and from other performed experiments that are not mentioned for the sake of brevity) it can be concluded that for $N>2$ both approximations have almost the same accuracy, our approximation being slightly better for small values of $\rho$, and the PCL approximation being slightly better for high values of $\rho$ (both methods are asymptotically exact as $\rho\uparrow1$). However, for $N=2$ our method suffers greatly from imbalance in the system, whereas the PCL approximation proves to be more robust. 
\begin{table}[h!t]
\parbox{0.45\textwidth}{
\begin{center}
\[
\begin{array}{|c|ccc|}
\hline
N  & 0-5\% & 5-10\% & 10-15\%  \\
\hline
2 &   89.32 & 9.11 & 1.56  \\
3 &    100.00 & 0.00 & 0.00  \\
4 &    100.00 & 0.00 & 0.00  \\
5 &    100.00 & 0.00 & 0.00  \\
\hline
\end{array}
\]
(a)
\end{center}
}
\hfill
\parbox{0.45\textwidth}{
\begin{center}
\[
\begin{array}{|c|ccc|}
\hline
N  & 0-5\% & 5-10\% & 10-15\%  \\
\hline
2 &  96.09 & 2.86 & 1.04  \\
3 &   99.31 & 0.69 & 0.00 \\
4 &   100.00 & 0.00 & 0.00 \\
5 &   100.00 & 0.00 & 0.00 \\
\hline
\end{array}
\]
(b)
\end{center}
}
\caption{Errors of the approximation applied to the 768 test cases with Poisson input, categorised in bins of $5\%$. In (a) the percentages of mean relative errors in each bin are shown for our approximation, in (b) results are shown for the PCL approximation.}
\label{numericalresults1poisson}
\end{table}

\subsection{Gated service}\label{gatedsubsection}

Until now we have only shown and discussed approximation results for polling systems with exhaustive service. The complete test bed described in Table \ref{testbed} has also been analysed for polling systems where each queue receives gated service. As can be seen in Table \ref{numericalresults1gated}, the overall quality of the approximation is good, but worse than for polling systems with exhaustive service. More details on the reason for these inaccuracies can be found in Table \ref{errorpervargated}, which is the equivalent of Table \ref{errorpervar} for gated service. Table \ref{errorpervargated}(b) illustrates that there is now a huge difference between systems with Poisson arrivals, and systems with non-Poisson arrivals. For the cases with $\C_{A_i}^2 = 1$, the approximation is extremely accurate, even for two-queue polling systems. The accuracy in cases with $\C_{A_i}^2 \neq 1$ is worse, which is caused by the assumptions that are made to approximate the LT limit \eqref{EWltgated}. Firstly, the decomposition \eqref{decompositiongated} does not hold for non-Poisson arrivals, and secondly, the terms $\E[I_i^\textit{res}]$ and $\frac{\E[V_iI_i]}{\E[I_i]}$ in this decomposition have only been approximated. For exhaustive service, these assumptions do not have much negative impact on the accuracy, but apparently, for gated service, they do. The mean relative errors for $N=2,\dots,5$ queues are respectively $2.70\%$, $2.25\%$, $1.90\%$, and $1.63\%$. The imbalance of the mean interarrival and service times hardly influences the accuracy of the approximation, as can be concluded from Table \ref{errorpervargated}(c).

If we consider the 768 cases with Poisson arrivals only, the mean relative errors of our approximation for $N=2,\dots,5$ are respectively $ 0.34\%, 0.17\%, 0.10\%$, and $0.08\%$. This accuracy is even better than the one achieved by the PCL approximation.
\begin{table}[h!t]
\[
\begin{array}{|c|cccc|}
\hline
N  & 0-5\% & 5-10\% & 10-15\% & 15-20\% \\
\hline
2 &   82.55 & 12.33 & 2.95 & 1.56 \\
3 &   85.42 & 10.53 & 3.13 & 0.81 \\
4 &   88.85 & 8.46 & 2.43 & 0.26 \\
5 &   92.22 & 6.60 & 1.15 & 0.03\\
\hline
\end{array}
\]
\caption{Errors of the approximation applied to the 2304 test cases with gated service, as described in Subsection \ref{gatedsubsection}, categorised in bins of $5\%$.}
\label{numericalresults1gated}
\end{table}

\begin{table}[h!t]
\begin{center}
\parbox{0.45\textwidth}{
\begin{center}
\[
\begin{array}{|c|cccccc|}
\hline
N & \multicolumn{6}{|c|}{\textrm{Load ($\rho$)}} \\
\hline
  & 0.10 & 0.30 & 0.50 & 0.70 & 0.90 & 0.99\\
\hline
2 &  2.64 & 4.55 & 4.31 & 3.10 & 1.25 & 0.37 \\
3 &  2.03 & 3.78 & 3.68 & 2.68 & 1.04 & 0.30 \\
4 &  1.62 & 3.14 & 3.13 & 2.32 & 0.92 & 0.28 \\
5 &  1.35 & 2.67 & 2.71 & 2.03 & 0.81 & 0.21 \\
\hline
\end{array}
\]
(a)
\end{center}
}
\hfill
\parbox{0.45\textwidth}{
\begin{center}
\[
\begin{array}{|c|ccccccccc|}
\hline
N & \multicolumn{9}{|c|}{\textrm{SCV interarrival times ($\C^2_{A_i}$)}} \\
\hline
  &&&  0.25 && 1    && 2    &&\\
\hline
2 &&&    4.72 && 0.34 && 3.05&& \\
3 &&&    4.06 && 0.17 && 2.53&& \\
4 &&&    3.45 && 0.10 && 2.16&& \\
5 &&&    2.98 && 0.08 && 1.84&&\\
\hline
\end{array}
\]
(b)
\end{center}
}
\[
\begin{array}{|c|cccc|}
\hline
N & \multicolumn{4}{|c|}{\textrm{Imbalance interarrival and service times}} \\
\hline
  & I_{A_i} = 1, I_{B_i} = 1 & I_{A_i} = 1, I_{B_i} = 5 & I_{A_i} = 5, I_{B_i} = 1 & I_{A_i} = 5, I_{B_i} = 5\\
\hline
2 &   2.76 & 2.64 & 2.81 & 2.59 \\
3 &   2.28 & 2.25 & 2.27 & 2.21 \\
4 &   1.93 & 1.91 & 1.90 & 1.87 \\
5 &   1.64 & 1.66 & 1.64 & 1.58 \\
\hline
\end{array}
\]
(c)
\end{center}
\caption{For gated service: mean relative approximation error, categorised by number of queues $(N)$ and total load of the system (a), SCV interarrival times (b), and imbalance of the interarrival and service times (c).}
\label{errorpervargated}
\end{table}

\section{Further research topics}\label{furtherresearch}

The research that is done in the present paper can be extended in many different directions. In this section we discuss some possibilities that we find most relevant.

\paragraph{Higher moments.}
Firstly, a logical follow-up step would be to use the same approach to find approximations for higher moments of the waiting time distribution as well. This might prove to be a hard exercise, since the LT limit of $\E[W_i^2]$ is unknown and, although its derivation might follow the same lines as in Section \ref{approximationsection}, it probably requires substantially more effort. In \cite{Kudoh96secondmoments}, explicit expressions for the second moments of the waiting time distributions are given, but only for symmetric systems with $N=2,3,$ and $4$, and under the assumption of Poisson arrivals. Also, the HT limit of $\E[W_i^2]$ is unknown, although some research in this area has already been done and in \cite{olsenvdmei05} a strong conjecture is given for the limiting distribution of $W_i$ as $\rho\uparrow 1$.

Another question that remains to be investigated, is the required form of the interpolation, as \eqref{ewapprox} is surely not adequate to approximate higher moments of $\E[W_i]$.

\paragraph{Other service disciplines.}
In the present paper, only exhaustive and gated service are discussed. In order to obtain results for polling systems with some queues receiving exhaustive service, and others receiving gated service, only minor modifications should be made, but we leave this to the reader. It would be more challenging to generalise the approximation to a wider variety of service disciplines. In particular, it would be nice to have one expression for the mean waiting time of customers in a queue with an arbitrary branching-type service discipline (cf. \cite{resing93}). The \emph{exhaustiveness} of a branching-type service discipline (cf. \cite{winandsPhD}) might appear in this expression. Gated and exhaustive are both branching type service disciplines, but are discussed separately in the present paper. The HT limit can most likely be established for arbitrary branching type service disciplines (see conjectures in \cite{olsenvdmei05}), so the question that remains is whether the LT limit can be found in a similar way.

\paragraph{Optimisation.}
One of the main reasons to choose \eqref{ewapprox} as form of the interpolation, besides its asymptotic correctness, is its simplicity. Having this exact and simple expression for the approximate mean waiting times, makes it very useful for optimisation purposes. In production environments, one can, for example, determine what the optimal strategy is to combine orders of different types (i.e., determine what queue customers should join). Because general arrivals are supported, one can determine optimal sizes of batches in which items are grouped and sent to a specific machine. The simplicity of \eqref{ewapprox} makes it possible for a manager to create a handy Excel sheet that can be used by operators to compute all kind of optimal parameter settings. No difficult computations are required at all, so a large variety of users can use the approximation.

In the present paper the accuracy of the approximation has been investigated and has been found to be very good in most situations. Another advantage of our approximation regarding optimisation purposes, is that the general shape of the approximated curve follows the exact curve very closely. Even in cases where the relative errors are rather large, like in Figure \ref{figExample1}, the shape of the actual curves is still very well approximated. This means that plugging our approximation, instead of an exact expression if it had been available, in an optimisation function yields an optimum that should be close to the true optimum.

\paragraph{Polling Table.}
The interpolation based approximation can also be extended to polling systems where the visiting order of the queues is not cyclic. Waiting times in polling systems with so-called polling tables can be obtained in the same way as shown in the present paper. Both the LT and HT limits are not difficult to determine in this situation, and the interpolation follows directly from these limits.

\paragraph{Model.}
The form of the interpolation might be changed to improve the accuracy of approximations for cases that give less satisfactory results in the present form. E.g., one could try other functions than a second-order polynomial as numerator of \eqref{ewapprox}. Alternatively, one could try to find a correction term which could be added to \eqref{ewapprox} to obtain better results for, e.g., two-queue polling systems. But most of all, if an \emph{exact} LT limit of the mean waiting time in a polling system with non-Poisson arrivals could be found, the accuracy of the approximation in the case of gated service might be improved.


\section*{Acknowledgements}

The authors wish to thank Onno Boxma for valuable discussions and for useful comments on earlier drafts of the present paper.

\bibliographystyle{abbrvnat}

\begin{thebibliography}{28}
\providecommand{\natexlab}[1]{#1}
\providecommand{\url}[1]{\texttt{#1}}
\expandafter\ifx\csname urlstyle\endcsname\relax
  \providecommand{\doi}[1]{doi: #1}\else
  \providecommand{\doi}{doi: \begingroup \urlstyle{rm}\Url}\fi

\bibitem[Bertsimas and Mourtzinou(1999)]{bertsimas99}
D.~Bertsimas and G.~Mourtzinou.
\newblock Decomposition results for general polling systems and their
  applications.
\newblock \emph{Queueing Systems}, 31:\penalty0 295--316, 1999.

\bibitem[Blanc and van~der Mei(1995)]{blancvdmei95}
J.~P.~C. Blanc and R.~D. van~der Mei.
\newblock Optimization of polling systems with {Bernoulli} schedules.
\newblock \emph{Performance Evaluation}, 22:\penalty0 139--158, 1995.

\bibitem[Boxma(1989)]{boxmaworkloadsandwaitingtimes89}
O.~J. Boxma.
\newblock Workloads and waiting times in single-server systems with multiple
  customer classes.
\newblock \emph{Queueing Systems}, 5:\penalty0 185--214, 1989.

\bibitem[Boxma and Groenendijk(1987)]{boxmagroenendijk87}
O.~J. Boxma and W.~P. Groenendijk.
\newblock Pseudo-conservation laws in cyclic-service systems.
\newblock \emph{Journal of Applied Probability}, 24\penalty0 (4):\penalty0
  949--964, 1987.

\bibitem[{Coffman, Jr.} et~al.(1995){Coffman, Jr.}, Puhalskii, and
  Reiman]{coffman95}
E.~G. {Coffman, Jr.}, A.~A. Puhalskii, and M.~I. Reiman.
\newblock Polling systems with zero switchover times: A heavy-traffic averaging
  principle.
\newblock \emph{The Annals of Applied Probability}, 5\penalty0 (3):\penalty0
  681--719, 1995.

\bibitem[{Coffman, Jr.} et~al.(1998){Coffman, Jr.}, Puhalskii, and
  Reiman]{coffman98}
E.~G. {Coffman, Jr.}, A.~A. Puhalskii, and M.~I. Reiman.
\newblock Polling systems in heavy-traffic: A {Bessel} process limit.
\newblock \emph{Mathematics of Operations Research}, 23:\penalty0 257--304,
  1998.

\bibitem[Everitt(1986)]{everitt86}
D.~Everitt.
\newblock Simple approximations for token rings.
\newblock \emph{IEEE Transactions on Communications}, COM-34\penalty0
  (7):\penalty0 719--721, 1986.

\bibitem[Fischer et~al.(2000)Fischer, Harris, and Xie]{fischer2000}
M.~J. Fischer, C.~M. Harris, and J.~Xie.
\newblock An interpolation approximation for expected wait in a time-limited
  polling system.
\newblock \emph{Computers \& Operations Research}, 27:\penalty0 353--366, 2000.

\bibitem[Fleming and Simon(1991)]{fleming1991}
P.~J. Fleming and B.~Simon.
\newblock Interpolation approximations of sojourn time distributions.
\newblock \emph{Operations Research}, 39\penalty0 (2):\penalty0 251--260, 1991.

\bibitem[Fuhrmann and Cooper(1985)]{fuhrmanncooper85}
S.~W. Fuhrmann and R.~B. Cooper.
\newblock Stochastic decompositions in the {$M/G/1$} queue with generalized
  vacations.
\newblock \emph{Operations Research}, 33\penalty0 (5):\penalty0 1117--1129,
  1985.

\bibitem[Groenendijk(1989)]{groenendijk89}
W.~P. Groenendijk.
\newblock Waiting-time approximations for cyclic-service systems with mixed
  service strategies.
\newblock In \emph{Proc. 12th ITC}, pages 1434--1441. North-Holland Publ. Co.,
  Amsterdam, 1989.

\bibitem[Keilson and Servi(1990)]{keilsonservi90}
J.~Keilson and L.~D. Servi.
\newblock The distributional form of {Little's Law} and the {Fuhrmann-Cooper}
  decomposition.
\newblock \emph{Operations Research Letters}, 9\penalty0 (4):\penalty0
  239--247, 1990.

\bibitem[Kudoh et~al.(1996)Kudoh, Takagi, and Hashida]{Kudoh96secondmoments}
S.~Kudoh, H.~Takagi, and O.~Hashida.
\newblock Second moments of the waiting time in symmetric polling systems.
\newblock \emph{Journal of the Operations Research Society of Japan},
  43\penalty0 (2):\penalty0 306--316, 1996.

\bibitem[Levy and Sidi(1990)]{levysidi90}
H.~Levy and M.~Sidi.
\newblock Polling systems: applications, modeling, and optimization.
\newblock \emph{IEEE Transactions on Communications}, 38:\penalty0 1750--1760,
  1990.

\bibitem[Olsen(2001)]{olsen2001}
T.~L. Olsen.
\newblock Limit theorems for polling models with increasing setups.
\newblock \emph{Probability in the Engineering and Informational Sciences},
  15\penalty0 (1):\penalty0 35--55, 2001.

\bibitem[Olsen and van~der Mei(2003)]{olsenvdmei03}
T.~L. Olsen and R.~D. van~der Mei.
\newblock Polling systems with periodic server routeing in heavy traffic:
  distribution of the delay.
\newblock \emph{Journal of Applied Probability}, 40:\penalty0 305--326, 2003.

\bibitem[Olsen and van~der Mei(2005)]{olsenvdmei05}
T.~L. Olsen and R.~D. van~der Mei.
\newblock Polling systems with periodic server routing in heavy traffic:
  renewal arrivals.
\newblock \emph{Operations Research Letters}, 33:\penalty0 17--25, 2005.

\bibitem[Reiman and Simon(1988)]{reimansimon88}
M.~I. Reiman and B.~Simon.
\newblock An interpolation approximation for queueing systems with {Poisson}
  input.
\newblock \emph{Operations Research}, 36\penalty0 (3):\penalty0 454--469, 1988.

\bibitem[Resing(1993)]{resing93}
J.~A.~C. Resing.
\newblock Polling systems and multitype branching processes.
\newblock \emph{Queueing Systems}, 13:\penalty0 409 -- 426, 1993.

\bibitem[Simon(1992)]{simon92}
B.~Simon.
\newblock A simple relationship between light and heavy traffic limits.
\newblock \emph{Operations Research}, 40\penalty0 (Supplement 2):\penalty0
  S342--S345, 1992.

\bibitem[Takagi(1988)]{takagi1988qap}
H.~Takagi.
\newblock Queuing analysis of polling models.
\newblock \emph{ACM Computing Surveys (CSUR)}, 20:\penalty0 5--28, 1988.

\bibitem[Tijms(1994)]{tijms94}
H.~C. Tijms.
\newblock \emph{Stochastic models: an algorithmic approach}.
\newblock Wiley, Chichester, 1994.

\bibitem[van~der Mei and Winands(2008)]{vdmeiwinands08}
R.~D. van~der Mei and E.~M.~M. Winands.
\newblock A note on polling models with renewal arrivals and nonzero
  switch-over times.
\newblock \emph{Operations Research Letters}, 36:\penalty0 500--505, 2008.

\bibitem[Vishnevskii and Semenova(2006)]{vishnevskiisemenova06}
V.~M. Vishnevskii and O.~V. Semenova.
\newblock Mathematical methods to study the polling systems.
\newblock \emph{Automation and Remote Control}, 67\penalty0 (2):\penalty0
  173--220, 2006.

\bibitem[Whitt(1989)]{whitt89}
W.~Whitt.
\newblock An interpolation approximation for the mean workload in a {$GI/G/1$}
  queue.
\newblock \emph{Operations Research}, 37\penalty0 (6):\penalty0 936--952, 1989.

\bibitem[Winands(2007{\natexlab{a}})]{winandsPhD}
E.~M.~M. Winands.
\newblock \emph{Polling, Production \& Priorities}.
\newblock PhD thesis, Eindhoven University of Technology, 2007{\natexlab{a}}.

\bibitem[Winands(2007{\natexlab{b}})]{winandslargesetups07}
E.~M.~M. Winands.
\newblock On polling systems with large setups.
\newblock \emph{Operations Research Letters}, 35:\penalty0 584--590,
  2007{\natexlab{b}}.

\bibitem[Winands(2009)]{winandslargesetupsbranching09}
E.~M.~M. Winands.
\newblock Branching-type polling systems with large setups.
\newblock \emph{To appear in OR Spectrum}, 2009.

\end{thebibliography}

\end{document}